\documentclass[12pt,a4paper]{amsart}

\usepackage{lscape}
\usepackage{adjustbox}
\usepackage{hyperref}

\newcommand {\fB}{{\mathfrak{B}}}
\newcommand {\fbgl}{{\mathfrak{bgl}}}

\newcommand {\fbrj}{{\mathfrak{brj}}}
\newcommand {\fel}{{\mathfrak{el}}}

\newcommand {\fwk}{{\mathfrak{wk}}}
\newcommand {\foo}{{\mathfrak{oo}}}

\newcommand {\DS}{{\text{DS}}}
\newcommand {\df}{{\text{df}}}
\newcommand {\ndf}{{\text{nds}}}

\DeclareMathOperator{\rank}{rank}

\usepackage{DLdef1}

\let\ssec\subsection

\renewcommand {\ssbegin}[2][*]
 {\refstepcounter{subsection}%
\if#1*
\addcontentsline{toc}{subsection}{\thesubsection.\hskip 1pc #2}%
\else
\addcontentsline{toc}{subsection}{\thesubsection.\hskip 1pc #2. #1}%
\fi
 \def \secno {\gdef \secno {}{\ssecfont
\thesubsection.\hskip 2ex}%
 }%
 \begin{#2}}

\renewcommand {\sssbegin}[2][*]
  {\refstepcounter{subsubsection}
\if#1*
\addcontentsline{toc}{subsubsection}{\thesubsubsection.\hskip 1pc #2}%
\else
\addcontentsline{toc}{subsubsection}{\thesubsubsection.\hskip 1pc #2. #1}
\fi
  \def \secno {\gdef \secno {}{\ssecfont \thesubsubsection.\hskip 2ex}%
  }%
   \begin{#2}}

\renewcommand {\parbegin}[2][*]
  {\refstepcounter{paragraph}
\if#1*
\addcontentsline{toc}{paragraph}{\theparagraph.\hskip 1pc #2}%
\else
\addcontentsline{toc}{paragraph}{\theparagraph.\hskip 1pc #2. #1}
\fi
  \def \secno {\gdef \secno {}{\ssecfont \theparagraph.\hskip 2ex}%
  }%
   \begin{#2}}

\makeatletter
\newcommand\Subsubsection[2]{%
\refstepcounter{subsubsection}%
\addcontentsline{toc}{subsubsection}{\thesubsubsection.\hskip 1pc #1. #2}%
\par\vskip1ex \@plus .4ex \@minus .2ex\noindent\thesubsubsection. \normalfont{\bfseries #1} (#2).\hskip.7em plus .3em}
\makeatother

\begin{document}

\title[Duflo--Serganova homology]{Duflo--Serganova homology for exceptional modular\\ Lie superalgebras with Cartan matrix}

\author{Andrey Krutov${}^{a,b,*}$, Dimitry Leites${}^{c, d}$, Jin Shang${}^c$}

\address{
 $^a$Department of Mathematics, University of Zagreb, Bijeni\v{c}ka 30, 10000 Zagreb, Croatia \\
${}^b$ Independent University of Moscow, Bolshoj Vlasievskij per, dom
 11, RU-119 002 Moscow, Russia;
 andrey.krutov@math.hr\\
${}^c$New York University Abu Dhabi,
Division of Science and Mathematics, P.O. Box 129188, United Arab
Emirates; $\{$dl146, js8544$\}$@nyu.edu\\
${}^d$Department of mathematics, Stockholm University, Roslagsv. 101, Stockholm, Sweden; mleites@math.su.se\\
${}^{*}$The corresponding author}


\begin{abstract} 
For the exceptional finite-dimensional modular Lie superalgebras $\mathfrak{g}(A)$ with indecomposable Cartan matrix $A$, and their simple subquotients, we computed non-isomorphic Lie superalgebras constituting the homologies of the odd elements with zero square.  These homologies are~key ingredients in the Duflo--Serganova approach to the representation theory.

There were two definitions of defect of Lie superalgebras in the literature with different ranges of application. We suggest a third definition and an easy-to-use way to find its value. 

In positive characteristic, we found out one more reason to consider the space of roots  over reals, unlike the space of weights, which should be considered over the ground field.

We proved that the  rank of the homological element (decisive in calculating the defect of a given Lie superalgebra)  should be considered in the adjoint module, not the irreducible module of least dimension (although the latter is sometimes possible to consider, e.g., for $p=0$). 

We also computed the above homology for the only case of simple Lie superalgebras with symmetric root system not considered so far over the field of complex numbers, and its  modular versions: $\mathfrak{psl}(a|a+pk)$ for $a$ and $k$ small, and $p=2, 3, 5$.

\end{abstract}


\maketitle

\markboth{\itshape 
A. Krutov\textup{,} D. Leites\textup{,} J. Shang}{{\itshape Duflo--Serganova homology}}

\thispagestyle{empty}
\setcounter{tocdepth}{3}
\tableofcontents

\section{Introduction}

For notation, see \cite{BGL, BGLd}. All vector spaces, in particular, Lie (super)algebras $\fg$ we consider here, are finite-dimensional 
over the ground field $\Kee$ of characteristic $p$. Let $p$ denote also the parity of the Lie superalgebra considered. In this paper, $\Kee^{m|n}$ denotes a~commutative Lie superalgebra, the one with zero bracket; $\fc$ denotes the center of $\fg$.

\ssec{Atypicality and defect} The irreducible finite-dimensional modules over the Lie superalgebra ${\fg=\fg(A)}$ with Cartan matrix $A$ over $\Cee$ are of the two types: generic  (Kac brightly called them ``typical"),  which are in one-to-one correspondence with certain $\fg_\ev$-modules, and ``atypical".  It soon became clear that different irreducible modules 
over simple  (and ``close'' to simple, like $\fgl$ is close to $\fpsl$) Lie superalgebras  have different levels of \textit{atypicality}, not exceeding the defined in \cite{KW} \textit{defect} $\df(\fg)$. Kac and Wakimoto defined $\df(\fg)$ only for Lie superalgebras of the form $\fg(A)$ over $\Cee$ with an \underline{\textbf{even}} non-degenerate invariant symmetric bilinear form (NIS for short).

The \textit{typical} irreducible $\fg(A)$-modules were first described analytically by Berezin in several cases (\cite{Ber}). Kac obtained a~complete classification of typical irreducible $\fg(A)$-modules using more adequate algebraic technique, see \cite{Ktyp} with corrections in \cite{Sg, SV}. 

The \textit{atypical} irreducible modules were first classified and described in terms of sections of vector bundles on certain Grassmann supervarieties of dimension $0|n$ by J.~Bernstein and D.~Leites for $\fsl(1|n)$ and $\fosp(2|2n)$ whose defect is equal to 1, see \cite{BL, L}. 

\ssec{The Duflo\/--Serganova functor} 
Over any ground field, let $\fg$ be a~Lie superalgebra, let
$M$ be a~$\fg$-module given by a~representation $\rho$. For an odd element $x\in\fg$ such that $x^2=0$ set
\be\label{g_x}
M_x := \Ker \rho_x/\IM \rho_x\text{~~and~~}\fg_x := \Ker \ad_x/\IM \ad_x.
\ee

It was known since long time ago that for any \textit{dg Lie superalgebra} with the odd differential~$d$, the space $M_d$ is a $\fg_d$-module.

This general fact was formulated usually for the non-zero differential $d$ of degree~$\pm 1$, see, e.g.,~\cite{Ge}. In these works, the Maurer-Cartan equation allows one to add only degree $\pm1$ elements $x\in\fg_\od$ to the differential~$d$.

Duflo and Serganova (see~\cite{DS}) considered a~completely different case where $d=0$, and the $\Zee$-grading of $\fg$ considered modulo 2 coincides with the parity of $\fg$. Then, one can add any odd element $x$ such that $x^2=0$ to the differential $d$,  and still have $(d+x)^2=0$.

Duflo and Serganova proved that the \textit{Duflo\/--Serganova functor}
$\DS_x : M \tto M_x$ from the category of $\fg$-modules to the category of $\fg_x$-modules is a~\textit{tensor} functor; this was a~new result. 

This functor helped to solve interesting problems,
see \cite{EASe, HPS, HR, IRS}, where $p=0$.

For the first applications of the DS-functor in the cases where $p>2$, see \cite{DK2}.

\textbf{For which $\fg$ was $\fg_x$ computed}? \underline{For $p=0$}, Duflo and Serganova (\cite{DS}) considered finite-dimensional Lie superalgebras $\fg$ of the form $\fg(A)$ with indecomposable and invertible Cartan matrix $A$, and $\fgl(a|b)$ for any $a, b$ over $\Cee$, i.e., all simple Lie superalgebras with symmetric root system, except $\fpsl(a|a)$. 

In this paper, we consider the cases of $\fpsl(a|a+kp)$ for $p=0, 2,3,5$, and $k$ and $a$ small.

The Lie superalgebras considered by Duflo and Serganova for $p=0$, as well as $\fpsl(a|a+kp)$ for $p>0$, and the exceptional simple Lie superalgebras and those of the form $\fg(A)$ we consider here have an even NIS, see \cite{BKLS, KLLS}.  

The specifically super analogs of $\fgl(n)$ --- Lie superalgebras $\fg$ of series $\fq$  their simple subquotiets,  have an \underline{\textbf{odd NIS} and no Cartan matrix}. For computation of the DS-homology of these Lie superalgebra, see \cite{KLS}. The case of Poisson Lie superalgebras $\fpo(0|n)$ and their simple subquotients $\fh'(0|n)$ will be considered separately.

\subsection{When $\fg_x\simeq\fg_y$ for $\fg=\fg(A)$ over $\Cee$} 
The bulk of \cite{DS} consists of the proof of the fact that, for the finite-dimensional Lie superalgebra $\fg=\fg(A)$ with indecomposable and \textbf{invertible} Cartan matrix $A$, the following statement holds:
\sssbegin[{\cite{DS}}: when $\fg_x\simeq\fg_y$]{Theorem}[{\cite{DS}}]\label{T1} We have
$\fg_{x}\simeq\fg_{y}\Longleftrightarrow\textup{rank}\ \rho_{x}=\textup{rank}\ \rho_{y}$, where
\be\label{these}
\begin{array}{l}
\rho=\begin{cases}\id\text{ in }V&\text{for $\fg$ of series $\fgl(V)$, or $\fsl(V)$ for $\sdim V\neq a|a$, or $\fosp(V)$}, \\
\ad\text{ in }\fg&\text{for $\fg=\fosp(4|2; a)$ or $\fag(2)$ or $\fab(3)$}.
\end{cases}
\end{array}
\ee
\end{Theorem}
In cases \eqref{these}, the module given by $\rho$ is the one of the least dimension, except for $\fosp(4|2; a)$, where $a=1$, or 2 or 3 (or the values of $a$ obtained from these under the action of $S_3$, see \cite{BGL}): for these values of $A$, the irreducible modules of the least superdimension are (up to the change of parity) $4|2$, $6|4$ and $8|6$, respectively. 

These exceptional modules of the least dimension are lucidly (but with a~typo) 
described in \cite{GL2}. For $a=2$ and 3, Duflo and Serganova proved Theorem~\ref{T1} without considering $\fosp(4|2; a)$-modules of the least dimension. We show that to prove that $\fg_{x}\simeq\fg_{y}$ we should compare ranks of the operators in the adjoint module, not in the ones of the least dimension. 

To show in which modules should one compute the rank of $\rho(x)$ in order to establish if $\fg_{x}\simeq\fg_{y}$, we do consider these modules of least dimension, see Table~\eqref{ospA}.

For the Lie superalgebra $\fg(A)$ with a~symmetrizable indecomposable Cartan matrix $A$, the space spanned by
roots $R$ over the ground field inherits the NIS $(-,-)$ given on $\fg(A)$, see \cite[eq.~(4)]{BKLS}. In the
modular case, however, we have to use a different definition of roots, most lucidly described in
\cite{BLLoS}. (Note that the NIS  on $\fg(A)$ induces a non-degenerate symmetric bilinear form on the space spanned by roots over $\Ree$ as
well.)

Observe that there is no NIS on $\fsl(n|n)$, but there is a NIS on $\fpsl(n|n)$ which has no $n|n$-dimensional irreducible modules; the $\fpsl(n|n)$-module of least dimensional is adjoint.

Recall that a~root is called \textit{odd} if the corresponding root vector is odd; $\beta\in R_\od$ is \textit{isotropic} if $(\beta, \beta) = 0$. Two roots $\alpha, \beta\in R$ are called \textit{orthogonal} if $(\alpha, \beta) = 0$. 

In addition to Theorem \ref{T1}, Duflo and Serganova proved the following 

\sssbegin[{\cite{DS}}: $\fg_x$ for $\fg=\fgl(m|n)$ and $\fsl(a|b)$, where $a\neq b$, and $\fosp(a|2b)$]{Theorem}[{\cite{DS}}: $\fg_x$ for $\fg=\fgl(m|n)$ and $\fsl(a|b)$, where $a\neq b$, and $\fosp(a|2b)$]\label{T2}
For  $\fg=\fgl(m|n)$ and $\fsl(a|b)$, where $a\neq b$, and $\fosp(a|2b)$ and for 
any set 
\be\label{df}
\{\beta_1,\dots,\beta_k\}\text{~~of linearly independent \textbf{mutually orthogonal} isotropic roots,}
\ee let
$x_k:=x_{\beta_1} +\dots+x_{\beta_k}$, and $\rho=\id$.

Then,
\be\label{x_k}
\begin{array}{l}
\rank \rho_{x_{k}}=\begin{cases}k&\text{for $\fgl(m|n)$ and $\fsl(m|n)$, where $m\neq n$},\\
2k&\text{for $\fosp(2m|2n)$ and $\fosp(2m+1|2n)$}, \end{cases}\\
\end{array}
\ee
and, respectively,
\be\label{g_x_k}
\begin{array}{l}
\fg_{x_k}\simeq \begin{cases}\text{$\fgl(m-k|n-k)$ and $\fsl(m-k|n-k)$},\\
\text{$\fosp(2m-2k|2n-2k)$ and $\fosp(2m-2k+1|2n-2k)$}. 
\end{cases}\\
\end{array}
\ee
\end{Theorem}

$\bullet$ For $p=0$, the Lie superalgebras $\fgl(n|n)$ is the only type of Lie superalgebras with non-invertible indecomposable Cartan matrix, see classification \cite{CCLL}. We consider the case of $\fpsl(n|n)$, the simple subquotient of $\fgl(n|n)$, in \S~\ref{Spsl}; the answer does not depend on $p$.

Kac and Wakimoto define $\df(\fg)$ as the maximal cardinality $k$ of the set~\eqref{df}. For $\fg=\fg(A)$ over $\Cee$, the defect $\df(\fg)$ is also equal to
\be\label{df1}
\text{the number $\ndf(\fg)$ of non-isomorphic DS-homology superalgebras $\fg_x$}.
\ee
Clearly, $\ndf(\fg)$ (same as $\text{def}$, see~ eq.~\eqref{def}) is independent of the existence of NIS on $\fg$.

Let us investigate: is there a~relation between the defect and the rank of $\rho_{x}$ in some irreducible representation $\rho$ apart from the one described in eq.~\eqref{x_k}? For example, for the simple Lie superalgebras $\fosp(4|2; a)$, the values of rank of $\rho_{x_{\ndf(\fg)}}$, where $\rho$ is either $\ad$ or the irreducible representation $\sigma$ of the least dimension, are as follows
\be\label{ospA}
\begin{tabular}{|l|c|c|c|}
\hline
$a$&rank$(\ad_{x_{\ndf(\fg)}})$&rank$(\sigma_{x_{\ndf(\fg)}})$&$\ndf(\fg)$\\
\hline
$1$&$8$&$2$&1 \\
\hline
$2$&$8$&$4$&1 \\
\hline
$3$&$8$&$6$&1 \\
\hline
\end{tabular}\ee

This table hints to define the defect looking at the adjoint representation instead of the one of the least dimension, but this hint is not as compelling, as the result of table~\eqref{!}. 

For the serial Lie superalgebras with Cartan matrix, both choices of $\rho$ (the tautological and adjoint) are OK, the tautological $\fg_x$-module in eq.~\eqref{g_x_k} is, however, of superdimension $\sdim V- \rank \rho_{x_{k}}|\rank \rho_{x_{k}}$, where $\rho$ is the representation in \textit{tautological} $\fg$-module $V$, not in the adjoint one.

For the series $\fpsl$, for $\fosp(4|2; a)$ for $a$ generic, and for the two exceptional simple Lie superalgebras --- $\fag(2)$ and $\fab(3)$ --- the irreducible representation of the least dimension is the adjoint one, so to define the rank of $x$ we have to take $\ad$ (and be happy we do not have to consider a~module of greater dimension).

\subsection{When $\fg_x\simeq\fg_y$ for $p>0$ and $\fg=\fg(A)$?}\label{ssQ} We want to answer the following questions:

$\bullet$ Q1) Should we take the space spanned by roots over the ground field $\Kee$, same as the space of weights, or should we follow the definition suggested by A.~Lebedev, see \cite{BLLoS}, and consider the inner product in the space spanned by the roots over $\Ree$, and not the one described in \cite{BKLS}?

$\bullet$ Q2) What should we take for $\rho$ trying to imitate the approach~\eqref{x_k}? (Does the answer to the question ``under what conditions on $x$ and $y$ we have $\fg_x\simeq\fg_y$?'' depend on $\rho$?)

The answers prove, by means of \textit{SuperLie} code, the following hypothesis that guided us.

\sssec{Hypothesis}\label{conj}
Among all Cartan matrices of $\fg(A)$ 
\be\label{Am}
\text{with the maximal number of $0$'s on its main diagonal}
\ee select the one, denote it $\cA$, whose Dynkin-Kac diagram $\cD$ has the maximal number $G_{max}$ of ``gray'' vertices not connected with each other (in terms of Cartan matrices this means that for all indices $i,j$ such that $\cA_{ii}=\cA_{jj}=0$ there is a maximal number of pairs $\cA_{ij}=\cA_{ji}=0$). 
Our Hypothesis consists of several parts:

1) \textbf{For any $\fg(A)$}, except for $\fbgl(4;a)$,  the non-isomorphic Lie (super)algebras $\fg_x$ correspond to $x=x_\beta$ for
$\beta=\beta_1+\dots+\beta_k$, where $k\leq G_{max}$ and each $\beta_i$ is one of the roots corresponding to the above ``grey''
vertices not connected with ``grey''
vertices corresponding to the other $\beta_j$. Moreover,\
\[
\df(\fg(A))=\ndf(\fg(A))=G_{max}.
\]

1e) For $\fg=\fbgl(4;a)$, we have $\ndf=G_{max}+1$ whereas $\df(\fg(A))=G_{max}$.

2)  \textbf{If $A$ is non-invertible}, the defect $\ndf(\fh)$ of the simple subquotient $\fh:=\fg(A)^{(1)}/\fc$ is given by one of the following formulas (\textbf{we were unable to find the pattern})
\[
\begin{array}{ll}
2a)& \text{$\ndf(\fh)=G_{max}-1$, unless $G_{max}=1$ in which case $\ndf(\fh)=G_{max}=1$,}\\
2b)&\text{$\ndf(\fh)=G_{max}$, e.g., $\ndf(\fpsl(n|n+pk))=\ndf(\fgl(n|n+pk))$ for $k\neq 0$.}\\
\end{array}
\]

Observe that for  $\fg=\fpsl(n|n)$, we have $\text{def}(\fg)=\ndf(\fg)\neq \df(\fg)$ in any characteristic $p$ we verified, see \S~3.

\Subsubsection{Fact}{When $\fg_x\simeq\fg_y$}\label{fact} For the \textbf{exceptional} simple Lie superalgebra $\fg$ of the form either $\fg(A)$, or $\fg(A)^{(1)}/\fc$, we have $\fg_x\simeq\fg_y$ if and only if $\rank \rho_x=\rank \rho_y$, where $\rho$ is the \textbf{adjoint}
representation of~$\fg$, \textbf{not} in the irreducible $\fg$-module $M$ of the least dimension, \emph{as
  examples show, e.g.,~ \eqref{bgl3;a},~\eqref{!}}.

\paragraph{Our Hypothesis is true over $\Cee$} Duflo and Serganova proved the statement \textit{equivalent to the claim of our Hypothesis} for Lie superalgebras $\fg(A)$ over $\Cee$ with Cartan matrix $A$, even non-invertible, as for $\fgl(a|a)$, see \cite{DS}. To see where the idea of the answer comes from, consider examples. 

\underline{Series $\fgl$ and $\fosp$}. To select the matrix $A_m$, see eq.~\eqref{Am}, we represent the supermatrices --- elements of $\fgl(a|b)$ --- in the \textit{alternating format} 
\[
\Par_{alt}=\begin{cases}(\ev,\od,\ev,\od,\dots, \ev,\od,\dots, \od)&\text{~~if $a\leq b$}\\
(\ev,\od,\ev,\od,\dots, \ev,\dots, \ev)&\text{~~if $a> b$}.
\end{cases} 
\]
Then, the elements $E_{i, i+1}$ (resp. $E_{i+1, i}$) are simple positive (resp. negative) root vectors for all~$i$. Since the neighboring root vectors do not commute, their roots are not orthogonal to each other; take every second root. 
So, $\ndf(\fgl(a|b))=\ndf(\fsl(a|b))=\min(a,b)$. 

The embeddings 
\[
\fgl(a|b)\tto\fosp(2a|2b)\tto\fosp(2a+1|2b)
\] 
give the value of $\ndf(\fg)$ for $\fg$ of the $\fosp$ series, see eq.~\eqref{x_k}: 
\[
\text{$\fgl(a|b)_{x_k}=\fgl(a-k|b-k)$ and $\fosp(N|2n)_{x_k}=\fosp(N-2k|2n-2k)$.}
\]

\underline{Exceptions: $\fg=\fosp(4|2; a)$ for $a$ generic, $\fag(2)$, and $\fab(3)$}. For them, the hypothesis describing when $\fg_x\simeq\fg_y$ is also true; in the hypothesis, $\rho$ is the adjoint representation, see eq.~\eqref{p>5}. For $\fg=\fosp(4|2; a)$ and $a=1,2$ or 3, we can take any of the two representations: either $\ad$, or the irreducible representation $\sigma$ of the least dimension, see eq.~\eqref{ospA}.

\subsection{Summary of our results}\label{OurRes}
We suggested a~simple method for 
determining the value of the defect of $\fg(A)$. We computed both $\ndf(\fg)$ and  \textit{the Duflo-Serganova homology}~$\fg_x$ in the two cases: 

(i)~over $\Kee$ of characteristic $p>0$, for each exceptional Lie superalgebra $\fg(A)$  with indecomposable Cartan matrix $A$; 

(ii)~for simple subquotients of the exceptional modular Lie superalgebras $\fg^{(1)}(A)/\fc$, where $A$ is not
invertible, e.g., for $\fpsl(n|n+pk)$ the conjecture is verified for $p=0, 2, 3, 5$, and $n$ and $k$  small.

Hypothesis~\ref{conj} is verified directly, by means of the \textit{SuperLie} package, see \cite{Gr}. Proving Hypothesis and answering Questions in Subsection~\ref{ssQ} we found out the following.

\Subsubsection{Fact}{Where do roots live, which representation determines defect}
\textup{F1)} Over $\Kee$, the space of roots should be considered over $\Ree$, as it was defined by A.~Lebedev, see \cite{BLLoS}, not over the ground field as the space of weights. For motivations, see Subsections $\ref{brj25}$, $\ref{g(1,6)}$.

\textup{F2)} To get the correct answer to the question ``when $\fg_x\simeq\fg_y$?'', we have to take $\rho=\ad$, \textbf{not} the irreducible representation $\sigma$ of the least dimension, see, e.g., table~\eqref{ospA}, and especially table~\eqref{!}: the same rank of $\sigma_x$ and $\sigma_y$ might correspond to $\fg_x\not\simeq\fg_y$, whereas $\text{rank}\ \sigma_x\neq\text{rank}\ \sigma_y$ might correspond to $\fg_x\simeq\fg_y$.

\section{The exceptional cases}\label{Sroot}

By $N\fg(A)$ we denote $\fg(A)$ corresponding to the $N$th Cartan matrix $A$ as listed in \cite{BGL}; recall that $\sdim A/a|B$ means that $\sdim \fg(A)= A|B$ and $\sdim \fg^{(1)}(A)/\fc=a|B$. The odd root vectors are \fbox{boxed} and isotropic roots are \underline{underlined}. The multiplication tables in $\fg_x$ are obtained with the aid of \textit{SuperLie} package. After several cases illustrating our answers, we refer the reader to \cite{BLLoS} for bulky lists of roots and the corresponding root vectors.

\ssec{Notation \protect$\fA\oplus_c \fB$, see \cite{BGL}}\label{4cases} This notation is needed to describe
the following Lie superalgebras or the corresponding DS-homologies 
\be\label{eq4cases}
\begin{array}{ll}
\text{$\fg(2,3)$, $\fg(2,6)$, and $\fg(3,3)$} & \text{for $p=3$},\\
\text{$\fbgl(4; \alpha)$,  $\fe(6, 6)$,  $\fe(7, 6)$, and $\fe(8, 1)$}&\text{for $p=2$}. \\
\end{array}
\ee
 This notation describes the case where $\fA$ and
$\fB$ are nontrivial central extensions of the Lie (super)algebras $\fa$
and $\fb$, respectively, and $\fA\oplus_c \fB$
--- a~nontrivial central extension of $\fa\oplus \fb$ (or, perhaps,
a more complicated semidirect sum $\fa\subplus \fb$, where $\fa$ is an ideal) with 1-dimensional center
spanned by $c$ --- is such that the restriction of the extension of
$\fa\oplus \fb$ to $\fa$ gives $\fA$ and that to $\fb$ gives $\fB$.

Consider the 4 Lie superalgebras $\fg(A)$, where $p=2$, listed in eq.~\eqref{eq4cases} in more details (see \cite{BGL}). Then, $\fg(A)_\ev$ is of the form
\be\label{GofB}
\fg(B)\oplus_c\fhei(2)\simeq \fg(B)\oplus\Span(X^+,X^-),
\ee 
where the
matrix $B$ is not invertible (so $\fg(B)$ has a~grading element $d$
and a~central element $c$), and where $X^+$, $X^-$ and $c$ span the
Heisenberg Lie algebra $\fhei(2)$. The brackets are:
\begin{equation}\label{strange}
\begin{array}{l}
{}[\fg^{(1)}(B), X^\pm]=0;\\
{}[d,X^\pm]=\begin{cases}X^\pm&\text{for $\fe(6, 6)$, $\fe(7,6)$, and $\fe(8,
1)$}; \\
\alpha X^\pm&\text{for $\fbgl(3; \alpha)$}\\
\end{cases}\\
{} [X^+,X^-]=c. \end{array}
\end{equation}

The odd part of $\fg(A)$ (at least in two of these four cases)
consists of two copies of the same $\fg(B)$-module $N$, the
operators $\ad_{X^\pm}$ permute these copies, and $\ad_{X^\pm}^2=0$,
so each of the operators maps one of the copies to the other, and
this other copy to zero.

\ssec{Notation \protect$\fA\oplus_c^d\fB$} Recall the definition and examples of \textit{double extensions}, see \cite{BLS}. Actually, in examples in eq.~\eqref{strange}, and~$\fg_x$ for~$\fe(7,6)$, we have \textit{double extensions}~$\fA$ and~$\fB$ of the Lie (super)algebras~$\fa$
and~$\fb$, respectively, and $\fA\oplus_c^d\fB$ is a double extension of~$\fa\oplus\fb$ by means of a central
element~$c$ and an outer derivation~$d$ such that the restriction of the extension of~$\fa\oplus\fb$ to~$\fa$ (resp.~$\fb$)
gives~$\fA$ (resp.~$\fB$). 

\ssec{For $p\geq 5$: $\fosp(4|2; a)$ for $a\neq 0, -1$, $\fag(2)$, and $\fab(3)$} The answer is the same as for $p=0$, namely: $\df(\fg)=1$ and $\fg_x$ is given by table~\eqref{p>5}; verified for $p=5, 7, 11$:
\begin{equation}\label{p>5}
\begin{tabular}{|l|l|l|l|}
\hline
$\fg$&$\fosp(4|2;a)$&$\fag(2)$&$\fab(3)$\\
\hline
$\fg_x$&$\Kee^{1|0}$&$\fsl(2)$ &$\fsl(3)$ \\
\hline
$\rank \ad_x$&$8$&$14$ &$16$ \\
\hline
\end{tabular}
\end{equation}

Each of the other exceptional Lie superalgebras $\fg(A)$ with indecomposable Cartan matrix $A$ exists only in characteristics 2, 3 and 5. The two Lie superalgebras $3\fg(2,3)$ and $1\fg(3,3)$ (indigenous to $p=3$) have ``the same''\footnote{Pretending that the elements of the Cartan matrix are integers, not elements of $\Kee$; we do the same  describing analogs of Serre relations, see \cite{BGL, BGLL}.} Cartan matrix as $3\fag(2)$ and $6\fab(3)$ (existing for $p=0$ and any $p>3$), respectively. Each of the other exceptional superalgebras $\fg(A)$ has no analogs among other exceptions, except for two pairs $\fbrj(2;5)\leftrightarrow\fbrj(2;3)$ and
$\fel(5;5)\leftrightarrow\fel(5;3)$, versions of which we consider one after the other for clarity. Although the algebras in these pairs have ``the same'' Cartan matrices for $p=5$ and 3, respectively, their structures are different (same thing with Lie superalgebras and their desuperizations for $p=2$, e.g., $\fbgl\leftrightarrow\fwk$). Let $\pi_i$ be the $i$th fundamental weight of the Lie algebra $\fg$, let $R(\sum \pi_i)$ denote the representation of $\fg$ with highest weight $\sum \pi_i$ and the corresponding $\fg$-module.

\sssec{$\fbrj(2;5)$ of $\sdim 10|12$}\label{brj25} We have $\fbrj(2; 5)_\ev = \fsp(4)$ and 
$\fbrj(2; 5)_\od = R(\pi_1 + \pi_2)$  is an irreducible the  $\fbrj(2; 5)_\ev$-module. We consider the following Cartan matrix and basis elements:
\be\label{brj25a}
\footnotesize
\begin{pmatrix}
0 & -1\\
-2 & 1\\
\end{pmatrix}\quad \tiny
\begin{tabular}{|l|l|} 
\hline
 the root vectors&the roots\\
\hline
\fbox{$x_{1}$},\;\fbox{$x_{2}$}&\underline{$\alpha_1$},\; $\alpha_2$\\
$x_3=\left[x_{1},\,x_{2}\right]$,\;
$x_4=\left[x_{2},\,x_{2}\right]$&$\alpha_1+\alpha_2,\ \ 2\alpha_2$\\
$\fbox{$x_5$}=\left[x_{2},\,\left[x_{1},\,x_{2}\right]\right]$&$\alpha_1+2\alpha_2$\\
$x_6=\left[\left[x_{1},\,x_{2}\right],\,\left[x_{2},\,x_{2}\right]
\right]$&$\alpha_1+3\alpha_2$\\
$\fbox{$x_7$}=\left[\left[x_{1},\,x_{2}\right],\left[x_{2},\,\left[x_{1},\,x_{2}\right]\right]\right],\;
\fbox{$x_8$}=\left[\left[x_{2},\,x_{2}\right],\left[x_{2},\,\left[x_{1},\,x_{2}\right]\right]\right]$&\underline{$2\alpha_1+3\alpha_2$},\ \ \underline{$\alpha_1+4\alpha_2$}\\
$x_9=\left[\left[x_{1},\,x_{2}\right],\left[\left[x_{1},\,x_{2}\right],\,\left[x_{2},\,x_{2}\right]
 \right]\right]$&$2\alpha_1+4\alpha_2$\\
$\fbox{$x_{10}$}=\left[\left[x_{2},\,\left[x_{1},\,x_{2}\right]
 \right],\,\left[\left[x_{1},\,x_{2}\right], \left[x_{2},\,x_{2}\right]\right]\right]$&\underline{$2\alpha_1+5\alpha_2$}\\
 \hline
\end{tabular}
\end{equation}

\be\label{brj25b}
\begin{array}{|l|l|c|l|}
\hline
x & \sdim\fg_x& \fg_x&\text{rank}\, \ad_x \\
\hline
x_1& 0|2 &\Kee^{0|2} &10\\
\hline
\end{array}
\ee
We get the same answer as for $x_1$ in eq.~\eqref{brj25a} for $x_7$, $x_8$, $x_{10}$, and $x_1+x_{10}$.

To find out which of the four homological linearly independent root vectors $x$'s in table~\eqref{brj25a} are orthogonal to one another (to compute the defect) we have to
decide do we consider the inner product in the space generated by the set of roots $R$ over the ground field $\Kee$ or over $\Ree$: unlike weights, which are considered over $\Kee$, the space of roots should be considered over $\Ree$ for several reasons as A.~Lebedev taught us, see \cite{BGL}; here is one more reason. 

Since $(\alpha_1, 2\alpha_1+5\alpha_2)\equiv 0\mod 5$, the only possible pair of mutually orthogonal isotropic candidates: $x_1$ and $x_{10}$; but since $\fg_x$ for $x=x_1+x_{10}$ is the same as for $x=x_1$, see table~\eqref{brj25a}, we conclude that $\ndf(\fg)=1$. Therefore, we have to 
\be\label{ALisRight}
\text{consider the inner product in the space of roots over $\Ree$.}
\ee
In this case, $(\alpha_1, 2\alpha_1+5\alpha_2) = - 10 \neq 0$. Thus $x_1$ and $x_{10}$ are not orthogonal.

\sssec{$1\fel(5;5)$ of $\sdim = 55|32$}\label{el(5;5)} We have $\fg_\ev = \fo(11)$ and $\fg_\od =R(\pi_5)$ as the $\fg_\ev$-module.
\be\label{el55a}
\begin{array}{|l|l|c|l|}
\hline
x & \sdim\fg_x & \fg_x&\text{rank}\,\ad_x \\
\hline
x_{3} &15 &\fsl(4)&32\\
\hline
\end{array}
\ee
We get the same answer as for $x_3$ in eq.~\eqref{el55a} for the other homological elements of the same rank
\[ 
x_{4},\
x_{5},\
x_{6},\
x_{7},\ 
x_{15},\ 
x_{17},\ 
x_{18},\ 
x_{19},\ 
x_{20},\ 
x_{26},\ 
x_{27},\ 
x_{28},\ 
x_{34},\ 
x_{35},\ 
x_{39}.
\]
Hypothesis~\ref{conj} is confirmed, $\ndf(\fel(5;3))=1$.

\ssec{For $p=3$}{}~{}
\sssec{$\fbrj(2;3)$ of $\sdim 10|8$}\label{ssbrj23} 
We have $\fg_\ev = \fsp(4)$ and 
$\fg_\od = R(2\pi_2)$ is an irreducible $\fg_\ev$-module. We consider the following Cartan matrix and basis elements
\be\label{brj23}
\begin{pmatrix}
0 & -1\\
-2 & 1\\
\end{pmatrix}\quad \tiny
\begin{tabular}{|l|l|} 
\hline
 the root vectors&the roots\\
\hline
\fbox{$x_{1}$},\;\fbox{$x_{2}$}&\underline{$\alpha_1$},\; $\alpha_2$\\
$x_3=\left[x_{1},\,x_{2}
\right]$, \ $x_4=\left[x_{2},\,x_{2}\right]$&$\alpha_1+\alpha_2$, \ $2\alpha_2$\\
$\fbox{$x_5$}=\left[x_{2},\,\left[x_{1},\,x_{2}
 \right]\right]$&$\alpha_1+2\alpha_2$\\
$x_6=\left[\left[x_{1},\,x_{2}\right],\left[x_{2},\,x_{2}\right]\right]$&$\alpha_1+3\alpha_2$\\
$\fbox{$x_7$}=\left[\left[x_{2},\,x_{2}\right],\left[x_{2},
\left[x_{1},\,x_{2}\right]\right]
 \right]$&\underline{$\alpha_1+4\alpha_2$}\\
$x_8=\left[
\left[x_{1},\,x_{2}\right],\left[\left[x_{1},\,x_{2}\right],
\left[x_{2},\,x_{2}\right]\right] \right]$&$2\alpha_1+4\alpha_2$\\
\hline
\end{tabular}
\end{equation}
\be\label{brj23a}
\begin{array}{|l|l|c|l|}
\hline
x & \sdim\fg_x & \fg_x&\text{rank}\,\ad_x \\
\hline
x_1, x_{7} & 2|0&\Kee^{2|0}&8\\
\hline
\end{array}
\ee

Since $(\alpha_1, \alpha_1+4\alpha_2)=-8$, it follows that $\ndf(\fg)=1$, Hypothesis~\ref{conj} is confirmed.

\sssec{$7\fel(5;3)$ of $\sdim = 39|32$}\label{el(5;3)}\label{el5;3} 
We have $\fg_\ev = \fo(9)\oplus \fsl(2)$ and $\fg_\od = R(\pi_4)\boxtimes \id$ as the $\fg_\ev$-module.

Let $M$ be the irreducible $\fg$-module with the highest weight $(0,0,0,0,1)$, then $\dim M = 18|16$; see~\cite{BGKL}. The answer in eq.~\eqref{el53} allows to consider $M$ as well as the adjoint module.
\be\label{el53}
\begin{array}{|l|l|c|l|l|}
\hline
x & \sdim\fg_x & \fg_x&\rank_\fg x & \rank_M x\\
\hline
x_{1} &15|8 &\fpsl(1|4) &24 & 10\\
\hline
x=x_4+x_5, x_1+x_3&7|0&\fpsl(3) & 32 & 16\\
\hline
\end{array}
\ee
We get the same answer as for $x_1$ in eq.~\eqref{el53} for the other homological elements of the same rank
\[ 
x_{2},\ 
x_{7},\ 
x_{10},\ 
x_{13},\ 
x_{15},\ 
x_{16},\ 
x_{20},\ 
x_{21},\ 
x_{24},\ 
x_{26},\ 
x_{27},\ 
x_{29},\ 
x_{30},\ 
x_{31},\ 
x_{33}.
\]
Hypothesis~\ref{conj} is confirmed, $\ndf(\fel(5;3))=2$.

\sssec{$1\fg(1,6)$ of $\sdim 21|14$}\label{g(1,6)} For $\fg=\fg(1,6)$, we have $\fg_\ev = \fsp(6)$ and $\fg_\od = R(\pi_3)$ as the $\fg_\ev$-module.
\be\label{g1,6a}
\begin{array}{|l|l|c|l|}
\hline
x & \sdim\fg_x& \fg_x&\rank_\fg x \\
\hline
x_3, \ x_9,\ x_{12},\ x_{15}& 7|0 &\fpsl(3)&14\\
\hline
\end{array}
\ee
Since $(\alpha_3, 2\alpha_1+4\alpha_2+\alpha_3)=-6\equiv 0\mod 3$, the only possible pair of mutually orthogonal isotropic candidates: $x_3$ and $x_{15}$; otherwise $\df(\fg)=1$. However, $\text{rank}\, \ad_{x_3+x_{15}}=\text{rank}\, \ad_{x_3}$. Therefore, we have to obey the rule~\eqref{ALisRight}. Hypothesis~\ref{conj} is confirmed.

\sssec{$2\fg(2,3)$ of $\sdim 12/10|14$}\label{sssG23} For $\fg=\fg(2,3)$, we have $\fg_\ev = \fgl(3) \oplus \fsl(2)$ and $\fg_\od = \fpsl(3)\boxtimes\id$ as the $\fg_\ev$-module. Clearly, $(\fg^{(1)}(2,3)/\fc)_\ev = \fpsl(3) \oplus \fsl(2)$.

\underline{For $\fg=\fg(2,3)$}, we have
\be
\begin{array}{|l|l|c|l|}
\hline
x & \sdim\fg_x & \fg_x&\rank_\fg x \\
\hline
x_{1} &2|4 & \fsl(1|1)\oplus_c\fsl(1|1)\oplus\Kee^{1|0} &10 \\
\hline
x_1+x_2& 1|3&\Kee^{1|3} &  11 \\
\hline
\end{array}
\ee
Hypothesis~\ref{conj} is confirmed, $\ndf(\fg)=2$.

\underline{For $\fh=\fg^{(1)}(2,3)/\fc$}, we have
\be
\begin{array}{|l|l|c|l|}
\hline
x & \sdim\fg_x & \fh_x&\rank_\fh x \\
\hline
x_{1} &0|4 &\Kee^{0|4} & 10 \\
\hline
x_1+x_2&0|4&\Kee^{0|4} & 10 \\
\hline
\end{array}
\ee
Hypothesis~\ref{conj} is confirmed, $\ndf(\fh)=1$.

\sssec{$2\fg(2,6)$ of $\sdim 36/34|20$}\label{g26} For $\fg=\fg(2,6)$, we have $\fg_\ev = \fgl(6)$ and $\fg_\od = R(\pi_3)$ as the $\fg_\ev$-module. Clearly, $(\fg^{(1)}(2,6)/\fc)_\ev = \fpsl(6)$.

$\fg = 2\fg(2,6)$
\be\label{g2,6a}
\begin{array}{|l|l|c|l|}
\hline
x & \sdim\fg_x& \fg_x&\rank_\fg x \\
\hline
x_2,\ x_3,\ x_4,\ & 16|0 & \fgl(3)\oplus_c^d\fgl(3) &20\\
\hline
\end{array}
\ee

$\fh= 2\fg^{(1)}(2,6)/\fc$
\be\label{g2,6a}
\begin{array}{|l|l|c|l|}
\hline
x & \sdim\fh_x& \fh_x&\rank_\fh x \\
\hline
x_2,\ x_3,\ x_4,\ & 14|0 &\fpsl(3)\oplus\fpsl(3) &20\\
\hline
\end{array}
\ee
We get the same answer for the other homological elements of the same rank. Hypothesis~\ref{conj} is confirmed, $\df(\fg)=\ndf(\fh)=1$.

\sssec{7$\fg(3,3)$ of $\sdim 23/21|16$}\label{SSSg33} 
Let $\spin_7:=R(\pi_3)$. We have 
\[
\text{$\fg(3,3)_\ev=(\fo(7)\supplus \Kee z)\subplus \Kee d$ and 
$\fg(3,3)_\od=(\spin_7)_+\oplus (\spin_7)_-$;}
\]
 the action of $d$ --- the outer derivative of $\fg(3,3)^{(1)}$ --- separates the identical $\fo(7)$-modules $\spin_7$ by acting on these modules as the scalar multiplication by $\pm1$, as indicated by subscripts, $z$ spans the center of $\fg(3,3)$. 
 
\underline{For $\fg = \fg(3,3)$}, we have 
\be\
\begin{array}{|l|l|c|l|}
\hline
x & \sdim\fg_x& \fg_x&\rank_\fg x \\
\hline
x_1,& 9|2 & \fgl(3)\oplus_c^d \fgl(1|1) &14\\
\hline
x_1+x_4& 7|0 &\fpsl(3) &16\\
\hline
\end{array}
\ee
Hypothesis~\ref{conj} is confirmed, $\ndf(\fg)=2$.

\underline{For $\fh = \fg^{(1)}(3,3)/\fc$}, we have 
\be\
\begin{array}{|l|l|c|l|}
\hline
x & \sdim\fh_x& \fh_x&\rank_\fh x \\
\hline
x_1,& 7|2 &\fpsl(3)\oplus \Kee^{0|2}&14\\
  \hline
x_1+x_4,& 7|2 &\fpsl(3)\oplus \Kee^{0|2}&14\\
\hline
\end{array}
\ee
Hypothesis~\ref{conj} is confirmed, $\ndf(\fh)=1$.

\sssec{$2\fg(3,6)$ of $\sdim 36|40$} For $\fg=\fg(3,6)$, we have $\fg_\ev = \fsp(8)$ and $\fg_\od = R(\pi_3)$ as the $\fg_\ev$-module.
\be\label{g3,6a}
\begin{array}{|l|l|c|l|}
\hline
x & \sdim\fg_x&\fg_x&\rank_\fg x \\
\hline
x_1 & 10|14 & \fg^{(1)}(2,3)/\fc& 26\\
\hline
x_1+x_4& 0|4 &\Kee^{0|4} &36\\
\hline
\end{array}
\ee
We get the same answer for the other homological elements of the same rank. Hypothesis~\ref{conj} is confirmed, $\ndf(\fg(3,6))=2$.

\sssec{$6\fg(4,3)$ of $\sdim 24|26$}\label{Lg4,3} For $\fg=\fg(4,3)$, we have $\fg_\ev = \fsp(6) \oplus \fsl(2)$ and $\fg_\od = R(\pi_2) \boxtimes \id$ as the $\fg_\ev$-module.
\be\label{g4,3a}
\begin{array}{|l|l|c|l|}
\hline
x & \sdim\fg_x& \fg_x&\rank_\fg x \\
\hline
x_1, x_2, x_3, x_{4}, x_{8} & 6|8 &\fpsl(2|2)&18\\
x_{12}, x_{17}, x_{20}, x_{22}& & &\\
\hline
x_1+x_3,\ x_2+x_4&0|2 &\Kee^{0|2}&24\\
\hline
\end{array}
\ee
We get the same answer for the other homological elements of the same rank. Hypothesis~\ref{conj} is confirmed, $\ndf(\fg)=2$.

\sssec{$13\fg(8,3)$ of $\sdim 55|50$}\label{e83} We have $\fg_\ev = \ff(4) \oplus \fsl(2)$ and $\fg_\od = R(\pi_4) \boxtimes \id$ as the $\fg_\ev$-module.
\be\label{g8,3a}
\begin{array}{|l|l|c|l|}
\hline
x & \sdim\fg_x&\fg_x&\rank_\fg x \\
\hline
x_2,\ x_3,\ x_{4},\ x_5 & 21|16 &\fg^{(1)}(2,3)/\fc&34\\

\hline
x_2+x_5,\ x_3+x_5& 7|2 &\fpsl(3)\oplus\Kee^{0|2} &48\\
\hline
\end{array}
\ee
We get the same answer for the other homological elements of the same rank. Hypothesis~\ref{conj} is confirmed, $\ndf(\fg)=2$.

\sssec{$2\fg(4,6)$ of $\sdim 66|32$}\label{g46} We have $\fg_\ev = \fo(12)$ and $\fg_\od = R(\pi_5)$ as the $\fg_\ev$-module.
\be\label{g4,6a}
\begin{array}{|l|l|c|l|}
\hline
x & \sdim\fg_x&\fg_x&\rank_\fg x \\
\hline
x_3& 34|0 &\fpsl(6)&32\\
\hline
\end{array}
\ee
We get the same answer for the other homological elements of the same rank. Hypothesis~\ref{conj} is confirmed, $\ndf(\fg)=1$.

\sssec{$4\fg(6,6)$ of $\sdim 78|64$}\label{g66} We have 
\[
\text{$\fg(6,6)_\ev= \fo(13)$ and $\fg(6,6)_\od = \spin_{13}:=R(\pi_6)$ as the $\fg_\ev$-module.}
\]
\be\label{g66a}
\begin{array}{|l|l|c|l|}
\hline
x & \sdim\fg_x& \fg_x&\text{rank}\, \ad_x \\
\hline
x_1& 34|20 &\fg^{(1)}(2,6)/\fc &34\\
\hline
x_1+x_{3}, x_1+x_{4} &14|0 &\fpsl(3)\oplus\fpsl(3)&48\\
\hline
\end{array}
\ee
We get the same answer for the other homological elements of the same rank. 
By Hypothesis~\ref{conj}, $\ndf(\fg)=2$, and we are done.

\sssec{$5\fg(8,6)$ of $\sdim 133|56$}\label{g86} We have $\fg_\ev = \fe(7)$ and $\fg_\od = R(\pi_1)$ as the $\fg_\ev$-module.
\be\label{g86}
\begin{array}{|l|l|c|l|}
\hline
x & \sdim\fg_x& \fg_x&\text{rank}\, \ad_x \\
\hline
x_2& 77|0 &\fe^{(1)}(6)/\fc &56\\
\hline
\end{array}
\ee
We get the same answer for the other homological elements of the same rank. Hypothesis~\ref{conj} is confirmed, $\ndf(\fg)=1$. (Recall that $\dim \fe(6)=79/77|0$ for $p= 3$, whereas $\dim \fe(6)=78$ for $p\neq 3$.)

\subsection{$p=2$} Numbering of vertices of the Dynkin diagram of $\fe$ type superalgebras goes along the string, starting from the end-vertex of the short branch, the end-vertex connected by an edge with the branching point is the last one.

\sssec{$\fbgl(3;\alpha)$, where $\alpha\neq 0, 1$; $\sdim=10/8|8$} The roots of $\fbgl(3;\alpha)$ are the same as those of $\fosp(4|2;\alpha)$ of $\sdim=9|8$ with the same division into even and odd ones.

We consider the following Cartan matrix and the corresponding
positive root vectors 
\be\label{bgl3}
\begin{pmatrix}
 0 & 1 & 0 \\
 1 & 0 & \alpha \\
 0 & \alpha &0
\end{pmatrix}
\quad \tiny
\begin{tabular}{|l|l|} 
\hline
 the root vectors&the roots\\
\hline
\fbox{$x_1$}, \ \fbox{$x_2$},\ \fbox{$x_3$},&$\underline{\alpha_1}$, \ $\underline{\alpha_2}$,\ $\underline{\alpha_3}$\\
$x_4=[x_1,\;x_2], \ x_5 =[x_2,\; x_3], $&$\alpha_1+\alpha_2, \ \alpha_2+\alpha_3,$\\
$\fbox{$x_6$}=[x_3, [x_1,\;x_2]]$&$\underline{\alpha_1+\alpha_2+\alpha_3}$\\
$x_7=[[x_1, x_2], [x_2, x_3]]$&$\alpha_1+2\alpha_2+\alpha_3$\\
\hline
\end{tabular}
\ee

In what follows, let $M$ be an irreducible highest weight $\fg$-module of the least dimension.

\underline{$\fg = \fbgl(3;a)$}, let $M$  be an irreducible module of $\sdim M = 4|4$, see \cite{BGKL}. We have
\be\label{bgl3a,etc}
\begin{tabular}{|l|l|l|l|}
 \hline
 $x$ & $\fg_x$ & $\rank_{\fg} x$& $\rank_M x$\\
 \hline
$x_1$ & $\Kee^{4|2}$ & $5$ & $2$\\
$x_2$ & $\Kee^{2|0}$ & $7$ & $3$\\
$x_3$ & $\Kee^{4|2}$ & $5$ & $2$\\
$x_1+x_3$ & $\Kee^{2|0}$ & $7$ & $4$\\ 
\hline\end{tabular}
\ee

\underline{For $\fh=\fbgl^{(1)}(3;\alpha)/\fc$}, the highest weight of~$M$ is~$(0,1,0)$ (with respect to the
maximal torus of $\fbgl(3;\alpha)$) and $\dim M = 8|6$, see \cite{BGKL}. 
\be\label{bgl3;a}
\begin{tabular}{|l|l|l|l|}
 \hline
 $x$ & $\fg_x$ & $\rank_{\fg} x$& $\rank_M x$\\
 \hline
 $x_1$ & $\Kee^{2|2}$ & $6$ & $8$ \\
 $x_2$ & $\Kee^{0|0}$ & $8$ & $6$ \\ 
 $x_3$ & $\Kee^{2|2}$ & $6$ & $4$ \\
 $x_1+x_3$ & $\Kee^{2|2}$ & $6$ & $6$ \\
 \hline
\end{tabular}
\ee

We see that $\rank_M x$ does not give a~conclusive information as to what $\fg_x$ is, unlike $\rank_\fg x$. Hypothesis~\ref{conj} item 1)  is confirmed: $\df(\fg)=2$, item 2) is NOT confirmed $\ndf(\fh)=2$.

\sssec{$\fbgl(4;\alpha)$, where $\alpha\neq 0, 1$;
$\sdim =18|16$} The roots of $\fbgl(4;\alpha)$ are the same as those of $\fwk(4;\alpha)$, but divided into even and odd ones. We consider the following Cartan matrix and the
corresponding positive root vectors
 \be\label{bgl4}\tiny
\begin{pmatrix}
 0 & \alpha & 0 & 0 \\
 \alpha & 0 & 1 & 0 \\
 0 & 1 & 0 & 1 \\
 0 & 0 & 1 & 0
\end{pmatrix}
\quad \tiny
\begin{tabular}{ll}
\fbox{$x_1$},\ \fbox{$x_2$} \ \fbox{$x_3$}, \ \fbox{$x_4$},&$\underline{\alpha_1}$, \ $\underline{\alpha_2}$, \ $\underline{\alpha_3}$, \ $\underline{\alpha_4}$\\
$x_5=[x_1,\;x_2]$,\ $x_6=[x_1, x_3]$, \ $x_7=[x_3,\;x_4],$&$\alpha_1+\alpha_2, \ \alpha_1+\alpha_3, \ \alpha_3+\alpha_4$\\
$\fbox{$x_8$}=[x_3,\;[x_1,x_2]],\ \fbox{$x_9$}=[x_4,\;[x_1,\;x_3]],$&$\underline{\alpha_1+\alpha_2+\alpha_3}, \ \underline{\alpha_1+\alpha_3+\alpha_4}$ \\
$x_{10}=[[x_1,x_2],\;[x_1,\;x_3]]$,\ $ x_{11}=[[x_1,\;x_2],[x_3,\;x_4]]$&$2\alpha_1+\alpha_2+\alpha_3, \ \alpha_1+\alpha_2+\alpha_3+\alpha_4 $\\
\fbox{$x_{12}$}=$[[x_1,\;x_2],\;[x_4,[x_1,\;x_3]]],$&$\underline{2\alpha_1+\alpha_2+\alpha_3+\alpha_4}
$\\
$x_{13}=[[x_3,\;[x_1,\;x_2]],[x_4,\;[x_1,\;x_3]]],$&$2\alpha_1+\alpha_2+2\alpha_3+\alpha_4$\\
$\fbox{$x_{14}$}=[[x_4,\;[x_1,x_3]],\;[[x_1,\;x_2],\;[x_1\text{,$\,
 $}x_3]]]$&$\underline{3\alpha_1+\alpha_2+2\alpha_3+\alpha_4} $\\
$x_{15}=[[[x_1,\;x_2],\;[x_1,\;x_3]],[[x_1,\;x_2],\;[x_3,x_4]]]$&$3\alpha_1+2\alpha_2+2\alpha_3+\alpha_4$
 \end{tabular}
\ee

Let $M$ be the irreducible $\fg$-module with highest weight $(0,0,0,1)$, then $\dim M = 16|16$, see \cite{BGKL}.
Let $V$ be the tautological $\fsl(3)$-module.

\be\label{!}
\begin{tabular}{|l|l|l|l|c|}
 \hline
 $x$ & $\dim \fg_x$ & $\rank_{\fg} x$& $\rank_M x$& $\fg_x$\\
 \hline
 $x_1$ & $8|6$ & $10$ & $10$ &$\fpsl(3|1)$\\
 $x_3$, $x_4$ & $8|6$ & $10$ & $8$ & $\fpsl(3|1)$\\
 $x_{14}$ & $8|6$ & $10$ & $16$ & $\fpsl(3|1)$\\
 $x_2$, $x_8$, $x_9$ & $4|2$ & $14$ & $12$
 & $\fg_x^{(3)} = 0$, $\fg_x/\fg_x^{(1)}=\Kee$, $\fg_x^{(1)}/\fg_x^{(2)}=\Kee^{2|2}$\\
 $x_{12}$ & $4|2$ & $14$ & $14$& $\fg_x^{(3)}=0$, $\fg_x/\fg_x^{(1)}=\Kee$, $\fg_x^{(1)}/\fg_x^{(2)}=\Kee^{2|2}$\\
 $x_1+x_3$ & $2|0$ & $16$ & $14$ & $\Kee^{2|0}$ \\
 $x_1+x_4$ & $2|0$ & $16$ & $14$ & $\Kee^{2|0}$ \\
 $x_2+x_4$ & $2|0$ & $16$ & $14$ & $\Kee^{2|0}$ \\
 \hline
\end{tabular}
\ee
Hypothesis~\ref{conj} is NOT confirmed: $\ndf(\fg)=3$. We see that $\rank_M x$ does not give a~conclusive information as to what $\fg_x$ is, unlike $\rank_\fg x$.

\sssec{$\fe(6,1)$ of $\sdim = 46|32$} We have
$\fg_\ev\simeq \fo\fc(2; 10)\oplus\Kee z$ and $\fg_\od$ is a
reducible module of the form $R(\pi_{4})\oplus R(\pi_{5})$ with the
two highest weight vectors: $y_5$ and 
\[
x_{36}=[[[x_4,x_5],[x_6,[x_2,x_3]]],
[[x_3,[x_1,x_2]],[x_6,[x_3,x_4]]]].
\]

Let $M$ be the irreducible $\fg$-module with highest weight $(1,0,0,0,0,0)$, then $\dim M = 16|11$, see \cite{BGKL}. The Cartan matrix of $\fe(6,1)$ we consider is that of $\fe(6)$ with the vector of parities of simple
roots being $111100$. 
\[
\begin{tabular}{|l|l|l|l|c|}
 \hline
 $x$ & $\dim \fg_x$ & $\rank_{\fg}x$& $\rank_M x$& $\fg_x$\\
 \hline

 $x_1,x_2,x_3,x_4,x_{10},x_{11},x_{12},x_{13},x_{17}$,
 & $24|10$ & $22$ & $6$ & $\fpsl(5|1)$ \\
 $x_{19},x_{20},x_{25},x_{27},x_{29},x_{31},x_{34}$
 & & & & \\
 \hline
 $x_1+x_3, \ x_2+x_4$ & $14|0$ & $32$ & $10$ & $\fpsl(4)$ \\
 \hline
\end{tabular}
\]
Hypothesis~\ref{conj} is confirmed, $\ndf(\fg)=2$. 

\sssec{$\fe(6,6)$ of $\sdim = 38|40$} 
In this case, $\fg(B)\simeq \fgl(6)$, see eq.~\eqref{GofB}. The module
$\fg_\od$ is irreducible with the highest weight vector
\[x_{35}=[[[x_3,x_6], [x_4,[x_2,x_3]]], [[x_4,x_5], [x_3,[x_1,x_2]]]]
\text{~~of weight $(0,0,1,0,0,1)$.}
\]

Let $M$ be the irreducible $\fg$-module with highest weight $(0,0,0,0,0,1)$, its dimension is $\dim M = 15|12$, see \cite{BGKL}. The Cartan matrix of $\fe(6,6)$ we consider is that of $\fe(6)$ with the vector of parities of simple
roots being $111111$. 
\[
\begin{tabular}{|l|l|l|l|c|}
 \hline
 $x$ & $\dim \fg_x$ & $\rank_{\fg} x$& $\rank_M x$& $\fg_x$\\
 \hline
 $x_1$ & $16|18$ & $22$ & $6$ & $\fpsl(3|3)$ \\
 \hline
 $x_1+x_3$ & $6|8$ & $32$ & $10$ & $\fpsl(2|2)$ \\
 \hline
 $x_1+x_3+x_5$ & $0|2$ & $38$ & $12$ & $\Kee^{0|2}$ \\
 \hline
\end{tabular}
\]
We get the same answer for the other elements of the same rank. Hypothesis~\ref{conj} is confirmed, $\ndf(\fg)=3$.

\sssec{$\fe(7,1)$ of $\sdim=80/78|54$} The Cartan matrix of $\fe(7,1)$ is that of $\fe(7)$ with the parities of simple roots $1111001$. 
\[
\begin{tabular}{|l|l|l|l|}
 \hline
 $x$           & $\dim \fg_x$ & $\rank_{\fg} x$& $\fg_x$\\
 \hline
 $x_1$         & $46|20$ & $34$ & $\foo^{(1)}_{\Pi\Pi}(2|10)$\\
 $x_1+x_3$     & $28|2$  & $52$ & $\fo^{(1)}_\Pi(8) \oplus_c^d \fgl(1|1)$\\
 $x_1+x_3+x_7$ & $26|0$  & $54$ & $\fo_\Pi^{(2)}(8)/\fc$\\
 \hline
\end{tabular}
\]

We get the same answer for the other elements of the same rank. Hypothesis~\ref{conj} is confirmed, $\ndf(\fg)=3$.

\underline{Consider $\fh:=(\fe^{(1)}(7,1))/\fc$}. We see that $\ndf(\fh)=2$, hypothesis~\ref{conj} is confirmed.
\[
\begin{tabular}{|l|l|l|l|}
 \hline
 $x$           & $\dim \fh_x$ & $\rank_{\fh} x$& $\fh_x$\\
 \hline
 $x_1$         & $44|20$ & $34$ & $\foo^{(2)}_{\Pi\Pi}(2|10)/\fc$\\
 $x_1+x_3$     & $26|2$  & $52$ & $\fo_\Pi^{(2)}(8)/\fc\oplus \Kee^{0|2}$\\
 $x_1+x_3+x_7$ & $26|2$  & $52$ & $\fo_\Pi^{(2)}(8)/\fc\oplus \Kee^{0|2}$\\
 \hline
\end{tabular}
\]

\sssec{$\fe(7,6)$ of $\sdim=70/68|64$}  The Cartan matrix of $\fe(7,6)$ we consider is that of $\fe(7)$ with the vector of parities of simple
roots being $0101010$. 
\[
\begin{tabular}{|l|l|l|l|}
 \hline
 $x$ & $\dim \fg_x$ & $\rank_{\fg} x$& $\fg_x$\\
 \hline
 $x_2$ & $36|30$ & $34$ & $\fpe^{(1)}(6)$\\ 
 $x_2+x_4$ & $18|12$ & $52$ & $\fpe^{(1)}(4)\oplus_{c}^{d} \fgl(2)$\\
 $x_2+x_4+x_6$ & $8|2$ & $62$ & solvable $\dim \fc = 1|2$,\\
 &&& $\dim \fg_x^{(i)}=\begin{cases}7|0&\text{if $i=1$,}\\
 1|0&\text{if $i=2$,}\\
 0|0&\text{if $i=3$,}\end{cases}$\\
 \hline
\end{tabular}
\]
We get the same answer for the other elements of the same rank. Hypothesis~\ref{conj} is confirmed, $\ndf(\fg)=3$.

\underline{Consider $\fh:=(\fe^{(1)}(7,6))/\fc$}. We see that $\ndf(\fh)=3$, hypothesis~\ref{conj} is confirmed.
\[
\begin{tabular}{|l|l|l|l|}
 \hline
 $x$ & $\dim \fh_x$ & $\rank_{\fh} x$& $\fh_x$\\
 \hline
 $x_2$ & $34|30$ & $34$ & $\fpe^{(2)}(6)/\fc$\\ 
 $x_2+x_4$ & $16|12$ & $52$ & $\fpe^{(2)}(4)/\fc\oplus \Kee^{2|0}$\\
 $x_2+x_4+x_6$ & $6|2$ & $62$ & $\Kee^{6|2}$\\
 \hline
\end{tabular}
\]

\sssec{$\fe(7,7)$ of $\sdim=64/62|70$}  The Cartan matrix of $\fe(7,7)$ we consider is that of $\fe(7)$ with the vector of parities of simple
roots being $1111111$.  
\[
\begin{tabular}{|l|l|l|l|l|}
 \hline
 $x$ & $\dim \fg_x$ & $\rank_{\fg} x$& $\fg_x$\\
 \hline
 $x_1$ & $30|36$ & $34$ & $\foo_{\Pi\Pi}^{(1)}(6|6)$\\
 $x_1+x_3$ & $12|18$ & $52$ & $\foo_{\Pi\Pi}^{(1)}(4|4)\oplus_c^d\fgl(1|1)$\\
 $x_1+x_3+x_5$ & $2|8$ & $62$ & solvable,\\
 &&& $\dim \fc=1|2$, $\dim \fg_x^{(i)}=\begin{cases}1|6&\text{if $i=1$,}\\
 1|0&\text{if $i=2$,}\\
 0|0&\text{if $i=3$,}\end{cases}$\\
 $x_1+x_3+x_5+x_7$ & $0|6$ & $64$ & $\Kee^{0|6}$ \\
 \hline
\end{tabular}
\]
We get the same answer for the other elements of the same rank. Hypothesis~\ref{conj} is confirmed, $\ndf(\fg)=3$.

\underline{Consider $\fh:=\fe^{(1)}(7,7)/\fc$}. We see that $\ndf(\fh)=3$, hypothesis~\ref{conj} is confirmed.

\[
\begin{tabular}{|l|l|l|l|l|}
 \hline
 $x$ & $\dim \fh_x$ & $\rank_{\fh} x$& $\fh_x$\\
 \hline
 $x_1$ &$28|36$ & $34$ & $\foo_{\Pi\Pi}^{(2)}(6|6)/\fc$\\
 $x_1+x_3$ & $10|18$ & $52$ & $\foo_{\Pi\Pi}^{(2)}(4|4)/\fc\oplus\Kee^{0|2}$\\
 $x_1+x_3+x_5$ & $0|8$ & $62$ & $\Kee^{0|8}$\\
 $x_1+x_3+x_5+x_7$ & $0|8$ & $62$ & $\Kee^{0|8}$\\
 \hline
\end{tabular}
\]

\sssec{$\fe(8,1)$ of $\sdim=136|112$} We have (see eq.~\eqref{GofB}) 
$\fg(B)\simeq \fe(7)$. The Cartan matrix of $\mathfrak{e}(8,1)$ we consider is that of $\fe(8)$ with the vector of parities of simple
roots being $11001111$. 
\[
\begin{tabular}{|l|l|l|l|l|}
 \hline
 $x$ & $\dim \fg_x$ & $\rank_{\fg} x$& $\fg_x$\\
 \hline
 $x_1$ & $78|54$ & $58$ & $\fe^{(1)}(7,1)/\fc$\\
 $x_1+x_5$ & $44|20$ & $92$ & $\foo^{(2)}_{\Pi\Pi}(2|10)/\fc$\\
 $x_1+x_5+x_7$ & $26|2$ & $110$ & $\fo^{(2)}_\Pi(8)/\fc\oplus\Kee^{0|2}$\\
$x_1+x_5+x_7+x_8$ &$6|2$ &$120$ & $\Kee^{6|2}$\\
 \hline
\end{tabular}
\]
We get the same answer for the other elements of the same rank. By Hypothesis~\ref{conj}, $\ndf(\fg)=4$.

\sssec{$\fe(8,8)$ of $\sdim=120|128$} In the $\Zee$-grading
with the Cartan matrix with the parities of simple roots
$00000001$ and with $\deg e_8^\pm=\pm 1$ and $\deg e_i^\pm=0$ for
$i\neq 8$, we have $\fg_0=\fgl(8)$. 

The Lie algebra $\fg_\ev$ is isomorphic to
$\fo_\Pi^{(2)}(16)\subplus\Kee d$, where
$d=E_{6,6}+\dots+E_{13,13}$, and $\fg_\od$ is an irreducible
$\fg_\ev$-module with the highest weight element
$x_{120}$ of weight $(1,0,\dots,0)$ with respect to $h_1,\dots,h_8$;
$\fg_\od$ also possesses a~lowest weight vector. The Cartan matrix of $\fe(8,8)$ we consider is that of $\fe(8)$ with the vector of parities of simple
roots being $11111111$. 
\[
\begin{tabular}{|l|l|l|l|l|}
 \hline
 $x$ & $\dim \fg_x$ & $\rank_{\fg} x$& $\fg_x$\\
 \hline
 $x_1$ & $62|70$ & $58$ & $\fe^{(1)}(7,7)/\fc$\\
 $x_1+x_3$ & $28|36$ & $92$ & $\foo_{\Pi\Pi}^{(2)}(6|6)/\fc$\\
 $x_1+x_3+x_5$ & $10|18$ & $110$ & $\foo_{\Pi\Pi}^{(2)}(4|4)/\fc\oplus\Kee^{0|2}$\\
$x_1+x_3+x_5+x_7$ & $0|8$ & $120$ & $\Kee^{0|8}$\\
 \hline
\end{tabular}
\]

We get the same answer for the other elements of the same rank. By Hypothesis~\ref{conj}, $\ndf(\fg)=4$.

\section{$\fpsl(n|n)$ for $p=0, 2,3,5$ (checked for $n=2,3, 4$)}\label{Spsl}
This case was not considered in \cite{DS,HR}. For it, the Hypothesis on the value of defect $\ndf$ is also true; the form of the answer differs from that in eq.~\eqref{g_x_k}. We consider the alternating format (so all Chevalley generators are odd).

For any $x\in \{x_1, x_3\}$, we have
\be\label{2,2}
\begin{array}{|l |l| l|}
\hline
\fg & \rank {\ad_x}, \ \fg_{x}&\rank \ad_{x_1 + x_3}, \ \fg_{x_1+x_3}\\
\hline
\fgl(2|2)&5, \ \fg_{x} =\fgl(1|1)
&7,\ \fg_{x_1+x_3} = 0\\
\hline
\fsl(2|2)&6,\ \fg_{x} = \fsl(1|1)\simeq\fhei(0|2)&
7, \ \fg_{x_1+x_3} = \Kee^{0|1}\\
\hline
\fpsl(2|2)&6, \ \fg_{x} = \Kee^{0|2}&
6, \ \fg_{x_1+x_3} = \Kee^{0|2}\\
\hline
\end{array}
\ee

For any $x\in \{x_1, x_3, x_{5}\}$, and $y=x_1+x_3$, and $z=x_1+x_3+x_5$, we have
\be\label{3,3}
\begin{array}{|l |l| l| l|}
\hline
\fg & \rank {\ad_x}, \ \fg_{x}&\rank \ad_{y}, \ \fg_{y}&\rank \ad_{z}, \ \fg_{z}\\
\hline
\fgl(3|3)&9, \ \fg_{x} =\fgl(2|2)
&15,\ \fg_{y} = \fgl(1|1)&17,\ \fg_{z} =0\\
\hline
\fsl(3|3)&10,\ \fg_{x} = \fsl(2|2)&
16, \ \fg_{y} = \fhei(0|2)&17,\ \fg_{z} = \Kee^{0|1}\\
\hline
\fpsl(3|3)&10, \ \fg_{x} =\fpsl(2|2)&
16, \ \fg_{y} = \Kee^{0|2}&16,\ \fg_{z} = \Kee^{0|2}\\
\hline
\end{array}
\ee

For any $x\in \{x_1, x_3, x_{5}, x_{7}\}$, and $y=x_1+x_3$, and $z=x_1+x_3+x_5$, and $v=x_1+x_3+x_5+x_7$, we have
\be\label{4,4}
\begin{array}{|l|l|l|l|l|}
\hline
\fg & \rank {\ad_x}, \ \fg_{x}&\rank \ad_{y}, \ \fg_{y}&\rank \ad_{z}, \ \fg_{z}&\rank \ad_{v},\ \fg_{v} \\
\hline
\fgl(4|4)&13, \ \fg_{x} =\fgl(3|3)
&23,\ \fg_{y} = \fgl(2|2)&29,\ \fg_{z} = \fgl(1|1)&31,\ \fg_{v}=0\\
\hline
\fsl(4|4)&14,\ \fg_{x} = \fsl(3|3)&
24, \ \fg_{y} = \fsl(2|2)&30,\ \fg_{z} = \fsl(1|1)&31,\ \fg_{v}=\Kee^{0|1}\\
\hline
\fpsl(4|4)&14, \ \fg_{x} =\fpsl(3|3)&
24, \ \fg_{y} = \fpsl(2|2)&30,\ \fg_{z} = \Kee^{0|2}&30,\ \fg_{v} = \Kee^{0|2}\\
\hline
\end{array}
\ee

\section{$\fpsl(n|n + pk)$ for $p=2,3,5$ and small non-zero $n, k$}\label{S4}

\subsection{$p=2$}
For $x = x_1$, and $y = x_1 + x_3$, we have
\be
\begin{array}{|l |l| l|}
 \hline 
 \fg & \rank {\ad_x}, \ \fg_{x}&\rank \ad_{y}, \ \fg_{y}\\
 \hline
  \fgl(2|4)  & 10,\ \ \fgl(1|3) & 16,\ \ \fgl(2)\\
 \fpsl(2|4) & 10,\ \ \fpsl(1|3)& 16,\ \ \Kee^{2|0}\simeq\fpsl(2)\\
 \hline
  \fgl(2|6)  & 14,\ \ \fgl(1|5) & 24,\ \ \fgl(5) \\
 \fpsl(2|6) & 14,\ \ \fpsl(1|5)& 24,\ \ \fpsl(4)\\
 \hline
  \fgl(2|8) &  18,\ \ \fgl(1|7) & 32,\ \ \fgl(6) \\ 
 \fpsl(2|8) & 18,\ \ \fpsl(1|7)& 32,\ \ \fpsl(6)\\
\hline
\end{array}
\ee

For $x = x_1$, $y = x_1 + x_3$, and $z=x_1+x_3+x_5$ we have
\be
\begin{array}{|l| l|l |l|}
 \hline 
 \fg & \rank {\ad_x}, \ \fg_{x}&\rank \ad_{y}, \ \fg_{y}&\rank \ad_{z},\ \ \fg_{z}\\
 \hline
  \fgl(3|5) & 14,\ \ \fgl(2|4) & 24,\ \ \fgl(1|3) & 30,\ \ \fgl(2)\\ 
 \fpsl(3|5)& 14,\ \ \fpsl(2|4)& 24,\ \ \fpsl(1|3)& 30,\ \ \Kee^{2|0}\simeq \fpsl(2)\\
  \hline
  \fgl(3|7) & 18,\ \ \fgl(2|6) & 32,\ \ \fgl(1|4) & 42,\ \ \fgl(4)\\ 
 \fpsl(3|7)& 18,\ \ \fpsl(2|6)& 32,\ \ \fpsl(1|5)& 42,\ \ \fpsl(4)\\
 \hline
  \fgl(3|9) & 22,\ \ \fgl(2|8) & 50,\ \ \fgl(1|7) & 54,\ \ \fgl(6)\\  
 \fpsl(3|9)& 22,\ \ \fpsl(2|8)& 40,\ \ \fpsl(1|7)& 54,\ \ \fpsl(6)\\
\hline
\end{array}
\ee

For $x = x_1$, $y = x_1 + x_3$, $z=x_1+x_3+x_5$, and $v=x_1+x_3+x_5+x_7$ we have
\be
\begin{array}{|l|l|l|l|l|}
 \hline 
 \fg & \rank {\ad_x}, \ \fg_{x}&\rank \ad_{y}, \ \fg_{y}&\rank \ad_{z},\ \ \fg_{z}
 & \rank \ad_{v},\ \ \fg_{v}\\
 \hline
  \fgl(4|6) & 18,\ \ \fgl(3|5) & 32,\ \ \fgl(2|4) & 42,\ \ \fgl(1|3) & 48,\ \ \fgl(2)\\
 \fpsl(4|6)& 18,\ \ \fpsl(3|5)& 32,\ \ \fpsl(2|4)& 42,\ \ \fpsl(1|3)& 48,\ \ \Kee^{2|0}\simeq\fpsl(2)\\
 \hline
  \fgl(4|8) & 22,\ \ \fgl(3|7) & 40,\ \ \fgl(2|6) & 54,\ \ \fgl(1|5) & 64,\ \ \fgl(4)\\ 
 \fpsl(4|8)& 22,\ \ \fpsl(3|7)& 40,\ \ \fpsl(2|6)& 54,\ \ \fpsl(1|5)& 64,\ \ \fpsl(4)\\
 \hline
  \fgl(4|10) & 26,\ \ \fgl(3|9) & 48,\ \ \fgl(2|8) & 66,\ \ \fgl(1|7) & 80,\ \ \fgl(6) \\ 
 \fpsl(4|10)& 26,\ \ \fpsl(3|9)& 48,\ \ \fpsl(2|8)& 66,\ \ \fpsl(1|7)& 80,\ \ \fpsl(6)\\
\hline
\end{array}
\ee

\subsection{$p=3$}

For $x = x_1$, and $y = x_1 + x_3$, we have
\be
\begin{array}{|l |l| l|}
 \hline 
 \fg & \rank {\ad_x}, \ \fg_{x}&\rank \ \ad_{y}, \ \fg_{y}\\
 \hline
  \fgl(2|5) & 12,\ \ \fgl(1|4) & 20,\ \ \fgl(3)\\
  \fpsl(2|5)& 12,\ \ \fpsl(1|4)& 20,\ \ \fpsl(3)\\
 \hline
  \fgl(2|8) & 18,\ \ \fgl(1|7) & 32,\ \ \fgl(6) \\
  \fpsl(2|8)& 18,\ \ \fpsl(1|7)& 32,\ \ \fpsl(6)\\
 \hline
  \fgl(2|11) & 24,\ \ \fgl(1|10) & 44,\ \ \fgl(9)\\ 
  \fpsl(2|11)& 24,\ \ \fpsl(1|10)& 44,\ \ \fpsl(9)\\
\hline
\end{array}
\ee

For $x = x_1$, $y = x_1 + x_3$, and $z=x_1+x_3+x_5$, we have
\be
\begin{array}{|l|l|l|l|}
 \hline 
 \fg & \rank {\ad_x}, \ \fg_{x}&\rank \ad_{y}, \ \fg_{y}&\rank \ad_{z},\ \ \fg_{z}\\
 \hline
 \fgl(3|6) & 16,\ \ \fgl(2|5)& 28,\ \ \fgl(1|4)& 36,\ \ \fgl(3)\\ 
 \fpsl(3|6) & 16,\ \ \fpsl(2|5)& 28,\ \ \fpsl(1|4)& 36,\ \ \fpsl(3)\\
 \hline
 \fgl(3|9) & 22,\ \ \fgl(2|8) & 40,\ \ \fgl(1|7) & 54,\ \ \fgl(6)\\
 \fpsl(3|9)& 22,\ \ \fpsl(2|8)& 40,\ \ \fpsl(1|7)& 54,\ \ \fpsl(6)\\
 \hline
 \fgl(3|12) & 28,\ \ \fgl(2|11) & 52,\ \ \fgl(1|10) & 72,\ \ \fgl(9)\\ 
 \fpsl(3|12)& 28,\ \ \fpsl(2|11)& 52,\ \ \fpsl(1|10)& 72,\ \ \fpsl(9)\\
\hline
\end{array}
\ee

For $x = x_1$, $y = x_1 + x_3$, $z=x_1+x_3+x_5$, and $v=x_1+x_3+x_5+x_7$, we have
\be
\begin{array}{|l|l|l|l|l|}
 \hline 
 \fg & \rank {\ad_x}, \ \fg_{x}&\rank \ad_{y}, \ \fg_{y}&\rank \ad_{z},\ \ \fg_{z}
 & \rank \ad_{v},\ \ \fg_{v}\\
 \hline
\fgl(4|7) & 20,\ \ \fgl(3|6) & 36,\ \ \fgl(2|5) & 48,\ \ \fgl(1|4) & 56,\ \ \fgl(3)\\
\fpsl(4|7)& 20,\ \ \fpsl(3|6)& 36,\ \ \fpsl(2|5)& 48,\ \ \fpsl(1|4)& 56,\ \ \fpsl(3)\\
 \hline
\fgl(4|10) & 26,\ \ \fgl(3|9) & 48,\ \ \fgl(2|8) & 66,\ \ \fgl(1|7) & 80,\ \ \fgl(6)\\ 
\fpsl(4|10)& 26,\ \ \fpsl(3|9)& 48,\ \ \fpsl(2|8)& 66,\ \ \fpsl(1|7)& 80,\ \ \fpsl(6)\\
\hline
\end{array}
\ee

\subsection{$p=5$}
For $x = x_1$, we have
\be
\begin{array}{|l| l|}
 \hline 
 \fg & \rank {\ad_x}, \ \fg_{x}\\
 \hline
 \fgl(1|6) & 12,\ \ \fgl(5) \\
 \fpsl(1|6)& 12,\ \ \fpsl(5)\\
 \hline
 \fgl(1|11) & 22,\ \ \fgl(10)\\ 
 \fpsl(1|11)& 22,\ \ \fpsl(10)\\
\hline
\end{array}
\ee

For $x = x_1$, and $y = x_1 + x_3$, we have
\be
\begin{array}{|l |l| l|}
 \hline 
 \fg & \rank {\ad_x}, \ \fg_{x}&\rank \ad_{y}, \ \fg_{y}\\
 \hline
 \fgl(2|7) & 16,\ \ \fgl(1|6) & 28,\ \ \fgl(5)\\
 \fpsl(2|7)& 16,\ \ \fpsl(1|6)& 28,\ \ \fpsl(5)\\
\hline
\fgl(2|12) & 26,\ \ \fgl(1|11) & 48,\ \ \fgl(10)\\ 
\fpsl(2|12)& 26,\ \ \fpsl(1|11)& 48,\ \ \fpsl(10)\\
\hline
\end{array}
\ee

For $x = x_1$, $y = x_1 + x_3$, and $z=x_1+x_3+x_5$, we have
\be
\begin{array}{|l|l|l|l|}
 \hline 
 \fg & \rank {\ad_x}, \ \fg_{x}&\rank \ad_{y}, \ \fg_{y}&\rank \ad_{z},\ \ \fg_{z}\\
 \hline
 \fgl(3|8) & 20,\ \ \fgl(2|7) & 36,\ \ \fgl(1|6) & 48,\ \ \fgl(5)\\
 \fpsl(3|8)& 20,\ \ \fpsl(2|7)& 36,\ \ \fpsl(1|6)& 48,\ \ \fpsl(5)\\
 \hline
\end{array}
\ee

For $x = x_1$, $y = x_1 + x_3$, $z=x_1+x_3+x_5$, and $v=x_1+x_3+x_5+x_7$, we have
\be
\begin{array}{|l|l|l|l|l|}
 \hline 
 \fg & \rank {\ad_x}, \ \fg_{x}&\rank \ad_{y}, \ \fg_{y}&\rank \ad_{z},\ \ \fg_{z}
 & \rank \ad_{v},\ \ \fg_{v}\\
 \hline
 \fgl(4|9) & 24,\ \ \fgl(3|8) & 44,\ \ \fgl(2|7) & 60,\ \ \fgl(1|6) & 72,\ \ \fgl(5)\\
 \fpsl(4|9)& 24,\ \ \fpsl(3|8)& 44,\ \ \fpsl(2|7)& 60,\ \ \fpsl(1|6)& 72,\ \ \fpsl(5)\\
\hline
\end{array}
\ee

\section{Comments}\label{Scomments}

\ssec{Two types of Lie superalgebras} The set of simple $\Zee$-graded Lie (super)algebras of finite dimension or of polynomial growth and their deformations is a~disjoint (at least, if $p>3$) union of two subsets: 
\begin{enumerate}
\item[(S)]\text{with a~symmetric set of roots relative the maximal torus,}\\
\text{i.e., with every root $\alpha$ there is a~root $-\alpha$ of the same multiplicity};\\[-4mm] 
\item[(N)]\text{with a~non-symmetric set of roots}.
\end{enumerate}

Some (or rather MOST) of the methods used to investigate Lie (super)algebras of type (S) rely on the existence of \textbf{Casimir elements} (in most cases, just one (degree-2) Casimir suffices) and symmetry of the root lattice. 

In the study of Lie (super)algebras of type (N) --- Lie (super)algebras of vector fields with polynomial or formal coefficients (briefly referred to as \textit{vectorial} Lie (super)algebras) and modules over them, these Casimirs are completely or partly absent\footnote{Such as the possibility to use the center of the universal enveloping algebra $U(\fg)$ which is trivial. However, in \cite{Ser}, Serganova showed that, at least for $\fg=\fpe(n)$, the superalgebra $U(\fg)$, whose center is trivial, should be replaced by $\overline U(\fg):=U(\fg)/\fr(U(\fg))$, whose center is sufficiently big, where $\fr(A)$ is the radical of the algebra $A$. It is very tempting to investigate applicability of Serganova's idea to other Lie (super)algebras $\fg$ with trivial center of $U(\fg)$. And with non-trivial centers as well.}, and hence can not be used. So, the problems of representation theory of the \lq\lq non-symmetric" (in particular, vectorial) Lie (super)algebras seem to be much more difficult than for the symmetric Lie (super)algebras.

Fortunately, over $\Cee$, the description of irreducible \textit{continuous} (with respect to the natural $(x)$-adic topology) modules over simple vectorial Lie algebras is very simple --- \textit{modulo representation theory of finite-dimensional simple Lie algebras}; superization is similar, see \cite{GLS}. 

\ssec{Two types of homological elements} Another tool for the study of Lie (super)algebras is  certain homology relative an element $x$. Here we considered an odd $x$ such that $x^2=2$.

For (N)-type simple modular Lie algebras for $p>2$, one considers the \textit{sandwich elements} $x$, i.e., such that $(\ad_x)^2=0$. S.~Kirillov proved that \textbf{the normalizer of the sandwich subalgebra}\footnote{In a~totally different setting and over $\Cee$, the term \textit{sandwich algebra} is used in \cite{Cu}, causing confusion.} \textbf{is the maximal subalgebra for $p>3$}, see \cite{Kir, KirS}. 

\sssec{On inhomogeneous ad-homological elements when $p=2$} We say that a non-zero $x\in \fg$ is \textit{ad-homological} if $(\ad_x)^2=0$. For $\fg=\fgl(2|6)$,
 inhomogeneous (with respect to parity) elements $x=x_1+x_3+x_5$ and $y=x_1+x_3+x_5+x_7$ are ad-homological. We get $\fg_x = \fgl(2)$ and $\fg_y = 0$, where $\rank_{\fg} \ad_x = 30$, and $\rank_{\fg} \ad_y = 32$. The meaning of ad-homological elements and their homology (when $p=2$) is unknown. Observe that  nobody computed the homology  corresponding to sandwiches; their meaning is also unknown.

Both sandwiches and DS-homology lead to what is called \textit{support varieties}. There are several non-equivalent definitions of these varieties: compare \cite{DS} with \cite{BaKN, BoKN, DK1, Ba, Ba1}. 

\ssec{Other definitions of the defect}\label{Bokn} \underline{In \cite{BoKN}}, there is given another definition of defect $\text{def}(\fg)$, equivalent to the above one $\df(\fg)$ for the simple (relatives of) Lie superalgebras $\fg$ with Cartan matrix over $\Cee$, but with a~wider range of application:
\be\label{def}
\text{def}(\fg):= \text{Krull dim}(H^{\bcdot}(\fg,\fg_\ev;\Cee))\text{~~for the trivial $\fg$-module $\Cee$}.
\ee

In \cite{BoKN}, it is mentioned that this definition is applicable to Lie superalgebras without Cartan matrix, such as $\fpe(n)$, $\fq(n)$ and their (not necessarily simple) relatives, and $\fgl(n|n)$. 

We observe that, moreover, the invariant~\eqref{def} is meaningful for \textbf{any} $\Zee/2$-grading of \textbf{any} Lie superalgebra and even for $\Zee/2$-graded Lie algebras. A priori, this cohomology is a~supercommutative superalgebra whose Krull dimension was recently defined, as on cue, see \cite{MZ}.

For $\fg=\fpe(n)$ and $\fq(n)$, and their relatives, $\text{def}(\fg)$ was computed in \cite{BoKN}, but $\fg_x$ was not computed for any $x$. For $\fq(n)$, and its simple subquotient, $\fg_x$ is computed in \cite{KLS}. 

\underline{The third definition of defect} --- $\ndf$ --- is given by formula \eqref{df1}.

\bigskip
\textbf{Acknowledgements}. For the possibility to conduct
difficult computations of this research we are grateful to M.~Al Barwani, Director of the High Performance
Computing resources at NYUAD. We thank E.~Herscovich for helpful comments. D.L. was partly supported by the grant AD 065 NYUAD.
A. K. was partly supported by the QuantiXLie Centre of Excellence, a~project cofinanced by the Croatian Government and European
Union through the European Regional Development Fund --- the Competitiveness and Cohesion Operational Programme (KK.01.1.1.01.0004).


\end{document}